\documentclass[10pt]{article}
\usepackage{amsfonts}
\usepackage[dvips]{graphicx}
\usepackage{amssymb,color}
\usepackage[latin1]{inputenc}
\usepackage{epic}
\usepackage[latin1]{inputenc}
\usepackage{enumerate}
\usepackage{color}
\usepackage{lettrine}
\usepackage{calc,pifont}
\usepackage[noindentafter,calcwidth]{titlesec}
\usepackage{amsmath}
\usepackage{amssymb}
\usepackage{amsthm}
\usepackage{subfigure}
\usepackage{cite}

\def\BBox{\kern  -0.2cm\hbox{\vrule width 0.2cm height 0.2cm}}

\newtheorem{lemma}{Lemma}[section]

\newtheorem{theorem}{Theorem}[section]

\newtheorem{corollary}{Corollary}[section]
\newtheorem{proposition}{Proposition}[section]

\newtheorem{example}[theorem]{Example}

\title{On transversal and 2-packing numbers in uniform linear systems}

\author{Carlos A. Alfaro \footnotemark[1] \and G. Araujo-Pardo \footnotemark[2] \and C. Rubio-Montiel \footnotemark[3] \and
	Adri{\' a}n V\'azquez-\'Avila  \footnotemark[4]}
\date{}

\begin{document}
\maketitle

\def\thefootnote{\fnsymbol{footnote}}
\footnotetext[1]{Banco de México, {\tt  carlos.alfaro@banxico.org.mx}.}

\footnotetext[2]{Instituto de Matemáticas, UNAM, {\tt garaujo@math.unam.mx}. This work was partially supported by PAPIIT-UNAM, Grant: IN107218 and  CONACyT-M\'exico, Grant: 282280}

\footnotetext[3]{Divisi{\' o}n de Matem{\' a}ticas e Ingenier{\' i}a, FES Acatl{\' a}n, UNAM, {\tt christian.rubio@apolo.acatlan.unam.mx}. This work was partially supported by PAIDI/007/19}

\footnotetext[4]{Subdirecci{\' o}n de Ingenier{\' i}a y Posgrado, UNAQ, {\tt adrian.vazquez@unaq.edu.mx}.}

\begin{abstract}
A linear system is a pair $(P,\mathcal{L})$ where
$\mathcal{L}$ is a family of subsets on a ground finite set $P$, such that $|l\cap l^\prime|\leq 1$, for every $l,l^\prime \in \mathcal{L}$. The elements of $P$ and $\mathcal{L}$ are called points and lines, respectively, and the linear system is called intersecting if any pair of lines intersect in exactly one point. A subset $T$ of points of $P$ is a transversal of $(P,\mathcal{L})$ if $T$ intersects any line, and the transversal number, $\tau(P,\mathcal{L})$, is the minimum order of a transversal. On the other hand, a 2-packing set of a linear system $(P,\mathcal{L})$ is a set $R$ of lines, such that any three of them have a common point, then the 2-packing number of $(P,\mathcal{L})$, $\nu_2(P,\mathcal{L})$, is the size of a maximum 2-packing set. It is known that the transversal number $\tau(P,\mathcal{L})$ is bounded above by a quadratic function of $\nu_2(P,\mathcal{L})$. An open problem is to haracterize the families of linear systems which satisfies $\tau(P,\mathcal{L})\leq \lambda\nu_2(P,\mathcal{L})$,  for some $\lambda\geq1$. In this paper, we give an infinite family of linear systems $(P,\mathcal{L})$ which satisfies $\tau(P,\mathcal{L})=\nu_2(P,\mathcal{L})$ with smallest possible cardinality of $\mathcal{L}$, as well as some properties of $r$-uniform intersecting linear systems $(P,\mathcal{L})$, such that $\tau(P,\mathcal{L})=\nu_2(P,\mathcal{L})=r$. Moreover, we state a characterization of $4$-uniform intersecting linear systems $(P,\mathcal{L})$ with $\tau(P,\mathcal{L})=\nu_2(P,\mathcal{L})=4$.
\end{abstract}




\section{Introduction}\label{sec:intro}

A \emph{linear system} is a pair $(P,\mathcal{L})$ where
$\mathcal{L}$ is a family of subsets on a ground finite set $P$, such that $|l\cap l^\prime|\leq 1$, for every pair of distinct subsets $l,l^\prime \in \mathcal{L}$. The linear system $(P,\mathcal{L})$ is \emph{intersecting} if $|l\cap l^\prime|=1$, for every pair of distinct subsets $l,l^\prime \in \mathcal{L}$. The elements of $P$ and $\mathcal{L}$ are called \emph{points} and \emph{lines}, respectively;  a line with exactly $r$ points is called a \emph{$r$-line}, and the \emph{rank} of $(P,\mathcal{L})$ is the maximum cardinality of a line in $(P,\mathcal{L})$, when all the lines of $(P,\mathcal{L})$ are $r$ lines we have a  \emph{$r$-uniform} linear system. In this context, a \emph{simple graph} is an 2-uniform linear system.

A subset $T\subseteq P$ is a \emph{transversal} (also called \emph{vertex cover} or \emph{hitting set} in many papers, as example \cite{MR1167472,Dorfling,MR2038482,MR1157424, 	MR3373359,MR2383447,MR2765536,MR1921545,MR1149871,MR1185788,Stersou}) of $(P,\mathcal{L})$ if for any line $l\in\mathcal{L}$ satisfies $T\cap l\neq\emptyset$. The \emph{transversal number} of $(P,\mathcal{L}),$ denoted by $\tau
(P,\mathcal{L}),$ is the smallest possible cardinality of a
transversal of $(P,\mathcal{L})$. 

A subset $R\subseteq \mathcal{L}$ is called \emph{2-packing} of $(P,\mathcal{L})$ if three elements are chosen in $R$ then they are not incident in a common point. The \emph{2-packing number} of $(P,\mathcal{L})$, denoted by
$\nu_2(P,\mathcal{L})$, is the maximum number of a
2-packing of $(P,\mathcal{L})$. 

There are many interesting works studying the relationship between these two parameters, for instance, in \cite{MR1185788}, the authors propose the problem of bounding $\tau(P,\mathcal{L})$ in terms of a function of $\nu_2(P,\mathcal{L})$ for any linear system. In \cite{MR3727901}, some authors of this paper and others proved that any linear system satisfies: 
\begin{equation}\label{desigualdad}
\lceil \nu_{2}/2\rceil\leq\tau
\leq \frac{\nu_2(\nu_2-1)}{2}.
\end{equation}
That is, the transversal number, $\tau$, of any linear system is upper bounded by a quadratic function of their 2-packing number, $\nu_2$.

In order to find how a function of $\nu_2(P,\mathcal{L})$ can bound $\tau(P,\mathcal{L})$, the authors of \cite{MR3021333} using probabilistic methods to prove that $\tau\leq\lambda\nu_2$ does not hold for any positive $\lambda$. In particular, they exhibit the existence of $k$-uniform linear systems $(P,\mathcal{L})$ for which their transversal number is $\tau(P,\mathcal{L})=n-o(n)$ and their $2$-packing number is upper bounded by $\frac{2n}{k}$. 

Nevertheless, there are some relevant works about families of linear systems in which their transversal numbers are upper bounded by a linear function of their 2-packing numbers. In \cite{CCA} the authors proved that if $(P,\mathcal{L})$ is a 2-uniform linear system, a simple graph, with $|\mathcal{L}|>\nu_2(P,\mathcal{L})$ then $\tau(P,\mathcal{L})\leq\nu_2(P,\mathcal{L})-1$; moreover, they characterize the simple connected graphs that attain this upper bound and the lower bound given in Equation (\ref{desigualdad}). In \cite{MR3727901} was proved that the linear systems $(P,\mathcal{L})$ with $|\mathcal{L}|>\nu_2(P,\mathcal{L})$ and $\nu_2(P,\mathcal{L})\in\{2, 3, 4\}$ satisfy $\tau(P,\mathcal{L})\leq\nu_2(P,\mathcal{L})$; and when attain the equality, they are a special family of linear subsystems of the projective plane of order $3$, $\Pi_3$, with transversal and $2$-packing numbers equal to $4$. Moreover, they proved that $\tau(\Pi_q)\leq\nu_2(\Pi_q)$ when $\Pi_q=(P_q,\mathcal{L}_q)$ is a projective plane of order $q$, consequently the equality holds when q is odd.  

The rest of this paper is structured as follows: In Section \ref{sec:desigualdad}, we present a result about linear systems satisfying $\tau\leq\nu_2-1$. In Section \ref{sec:familia_igualdad}, we give an infinite family of linear systems such that $\tau=\nu_2$ with smallest possible cardinality of lines. And, finally, in the last section, we presented some properties of the $r$-uniform linear systems, such that $\tau=\nu_2=r$, and we characterize the $4$-uniform linear systems with $\tau=\nu_2=4$.

\section{On linear systems with $\tau\leq\nu_2-1$}\label{sec:desigualdad}

Let $(P,\mathcal{L})$ be a linear system and $p\in P$ be a point. It is denoted by $\mathcal{L}_p$ to the set of lines incident to $p$. The \emph{degree} of $p$ is defined as $deg(p)=|\mathcal{L}_p|$ and the maximum degree overall points of the linear systems is denoted by
$\Delta(P,\mathcal{L})$. A point
of degrees $2$ and $3$ is called \emph{double} and \emph{triple} point, respectively, and two points $p$ and $q$ in $(P,\mathcal{L})$ are \emph{adjacent} if there is a line $l\in\mathcal{L}$ with $\{p,q\}\subseteq l$.

In this section, we generalize Proposition 2.1, Proposition 2.2, Lemma 2.1, Lemma 3.1 and Lemma 4.1 of \cite{MR3727901} proving that a linear system $(P,\mathcal{L})$ with $|\mathcal{L}|>\nu_2(P,\mathcal{L})$ and ``few'' lines satisfies $\tau(P,\mathcal{L})\leq\nu_2(P,\mathcal{L})-1$. Notice that, through this paper, all linear systems $(P,\mathcal{L})$ are considered with $|\mathcal{L}|>\nu_2(P,\mathcal{L})$ due to the fact $|\mathcal{L}|=\nu_2(P,\mathcal{L})$ if and only if $\Delta(P,\mathcal{L})\leq2$.

\begin{theorem}\label{thm:casicucarachageneral}
	Let $(P,\mathcal{L})$ be a linear system with $p,q\in P$ be two points such that $deg(p)=\Delta(P,\mathcal{L})$ and $deg(q)=\max\{deg(x): x\in
	P\setminus\{p\}\}$. If $|\mathcal{L}|\leq deg(p)+deg(q)+\nu_2(P,\mathcal{L})-3$, then $\tau(P,\mathcal{L})\leq\nu_2(P,\mathcal{L})-1$.
\end{theorem}
\paragraph{Proof}
Let $p,q\in P$ be two points as in the theorem, and let $\mathcal{L}''=\mathcal{L}\setminus\{\mathcal{L}_p\cup\mathcal{L}_q\}$, which implies that 
$|\mathcal{L}''|\leq\nu_2(P,\mathcal{L})-2$. Assume that $|\mathcal{L}''|=\nu_2(P,\mathcal{L})-2$ ($\mathcal{L}_p\cap \mathcal{L}_q\neq\emptyset$), otherwise,
the following set $\{p,q\}\cup\{a_l:a_l\mbox{ is any point of }l\in\mathcal{L}''\}$ is a transversal of $(P,\mathcal{L})$ of cardinality at most $\nu_2(P,\mathcal{L})-1$, and the statement holds. Suppose that $\mathcal{L}''=\{L_1,\ldots,L_{\nu_2-2}\}$ is a set of pairwise disjoint lines because, in other\-wise, they induce at least a double point, $x\in P$, hence the following set of points $\{p,q,x\}\cup\{a_l:l\in\mathcal{L}''\setminus \mathcal{L}''_x\}$, where $a_l$ is any point of $l$, is a transversal of $(P,\mathcal{L})$ of cardinality at most $\nu_2(P,\mathcal{L})-1$, and the statement holds.	
Let $l_q\in\mathcal{L}_q\setminus\{l_{p,q}\}$ be a fixed line and let $l_p$ be any line of $\mathcal{L}_p\setminus\{l_{p,q}\}$, where $l_{p,q}$ is the line containing to $p$ and $q$ (since $\mathcal{L}_p\cap \mathcal{L}_q\neq\emptyset$). Then $l_p\cap l_q\neq\emptyset$, since the $l_q$ induce a triple point on the following 2-packing $\mathcal{L}''\cup\{l_p,l_{p,q}\}$, which implies that there exists a line $L_{p,q}\in\mathcal{L}''$ with $l_q\cap l_p\cap L_{p,q}\neq\emptyset$, and hence $l_p\cap l_q\neq\emptyset$. Consequently, $deg(q)=\Delta(P,\mathcal{L})$ and $\Delta(P,\mathcal{L})\leq\nu_2(P,\mathcal{L})-1$ (since $\deg(p)-1\leq\nu_2(P,\mathcal{L})-2$). Therefore, the following set: $$\{l_p\cap L_i:i=1,\ldots,\Delta-1\}\cup\{a_\Delta,\ldots,a_{\nu_2-2}\}\cup\{p\},$$where $a_i$ is any point of $L_i$, for $i=\Delta,\ldots,\nu_2-2$, is a transversal of $(P,\mathcal{L})$ of the cardinality at most $\nu_2(P,\mathcal{L})-1$, and the statement holds.\qed
\section{A family of uniform linear systems with $\tau=\nu_2$}\label{sec:familia_igualdad}

In this section, we exhibit an infinite family of linear systems $(P,\mathcal{L})$ with two points of maximum degree and $|\mathcal{L}|=2\Delta(P,\mathcal{L})+\nu_2(P,\mathcal{L})-2$ with $\tau(P,\mathcal{L})=\nu_2(P,\mathcal{L})$. It is immediately, by Theorem \ref{thm:casicucarachageneral}, that  $\tau(P,\mathcal{L})\leq\nu_2(P,\mathcal{L})-1$ for linear systems with less lines.

In the remainder of this paper, $(\Gamma,+)$ is an additive Abelian group with neutral element $e$. Moreover, if $\sum_{g\in\Gamma}g=e$, then the group is called \emph{neutral sum group}. In the following, every group $(\Gamma,+)$ is a neutral sum group, such that $2g\neq e$, for all $g\in\Gamma\setminus\{e\}$. As an example of this type of groups we have $(\mathbb{Z}_n,+)$, for $n\geq3$ odd.

Let $n=2k+1$, with $k$ a positive integer, and $(\Gamma,+)$ be a neutral sum group of order $n$. Let:
$$\mathcal{L}=\{L_g:g\in\Gamma\setminus\{e\}\},\mbox{ where } L_g=\{(h,g):h\in\Gamma\},$$for $g\in\Gamma\setminus\{e\}$, and: $$\mathcal{L}_{p}=\{l_{p_g}:g\in\Gamma\},\mbox{ where } l_{p_g}=\{(g,h):h\in\Gamma\setminus\{e\}\}\cup\{p\},$$ for $g\in\Gamma$, and $\mathcal{L}_q=\{l_{q_g}:g\in\Gamma\}$, where:
$$l_{q_g}=\{(h,f_g(h)):h\in
\Gamma, f_g(h)=h+g \mbox{ with $f_g(h)\neq e$}\}\cup\{q\},$$ for $g\in\Gamma$.

Hence, the set of lines $\mathcal{L}$ is a set of pairwise disjoint lines with $|\mathcal{L}|=n-1$ and each line of $\mathcal{L}$ has $n$ points. On the other hand, $\mathcal{L}_p$ and $\mathcal{L}_q$ are set of lines incidents to $p$ and $q$, respectively, with $|\mathcal{L}_p|=|\mathcal{L}_p|=n$, and each line of $\mathcal{L}_p\cup\mathcal{L}_q$ has $n$ points. Moreover, this set of lines satisfies that, giving $l_{p_a}\in\mathcal{L}_p$ there exists an unique $l_{q_b}\in \mathcal{L}_q$ with $l_{p_a}\cap l_{q_b}=\emptyset$, otherwise, there exits $l_{p_a}\in \mathcal{L}_p$ such that
$l_{p_a}\cap l_{q_b}\neq\emptyset$, for all $l_{q_b}\in \mathcal{L}_q$, which implies that $a+b\in\Gamma\setminus\{e\}$, for all $b\in\Gamma$, which is a contradiction.

The linear system $(P_n,\mathcal{L}_n)$ with
$P_n=(\Gamma\times\Gamma\setminus\{e\})\cup\{p,q\}\quad\!\!\!\!\mbox{and}\quad\!\!\!\!\mathcal{L}_n=\mathcal{L}\cup\mathcal{L}_p\cup\mathcal{L}_q$, denoted by $\mathcal{C}_{n,n+1}$, is an  $n$-uniform linear system with $n(n-1)+2$ points and $3n-1$ lines. Notice that, this linear system has 2 points of degree $n$ (points $p$ and $q$) and $n(n-1)$ points of degree $3$.

A {\emph{linear subsystem}} $(P^{\prime },\mathcal{L}^{\prime })$ of a linear system $(P,\mathcal{L})$ satisfies that for any line $l^\prime\in \mathcal{L}^\prime$ there exists a line $l\in\mathcal{L}$ such that $l^\prime=l\cap P^\prime$, where $P^\prime\subset P$. Given a linear system $(P,\mathcal{L})$ and a point $p\in P$, the linear system obtained from $(P,\mathcal{L})$ by \emph{deleting the point $p$} is the linear system $(P^{\prime },\mathcal{L}^{\prime })$ induced by $\mathcal{L}^{\prime }=\{l\setminus \{p\}: l\in \mathcal{L}\}$. On the other hand, given a linear system $(P,\mathcal{L})$ and a line $l\in \mathcal{L}$, the linear system obtained from $(P,\mathcal{L})$ by \emph{deleting the line $l$} is the linear system $(P^{\prime },\mathcal{L}^{\prime })$ induced by $\mathcal{L}^{\prime }=
\mathcal{L}\setminus \{l\}$. The linear systems $(P,\mathcal{L})$ and $(Q,\mathcal{M})$ are isomorphic, denoted by $(P,\mathcal{L})\simeq(Q,\mathcal{M})$, if after deleting the points of degree 1 or 0 from both, the systems $(P,\mathcal{L})$ and $(Q,\mathcal{M})$ are isomorphic as hypergraphs (see \cite{berge1984hypergraphs}).

It is important to state that in the rest of this paper it is considered linear systems $(P,\mathcal{L})$ without points of degree one because, if $(P,\mathcal{L})$ is a linear
system which has all lines with at least two points of degree 2 or more, and $(P^\prime,\mathcal{L}^\prime)$ is the linear system obtained from $(P,\mathcal{L})$ by deleting all points of degree one, then they are essentially the same linear system because it is not difficult to prove that transversal and 2-packing numbers of both coincide (see \cite{MR3727901}).
\begin{figure}
	\begin{center}
		\includegraphics[height =4cm]{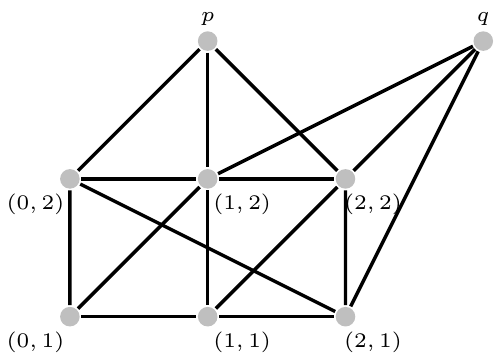}
	\end{center}
	\caption{Linear system $\mathcal{C}_{3,4}=(P_3,\mathcal{L}_3)$.}\label{fig:nu_2_3}
\end{figure}   
\begin{example}
	Let $\Gamma=\mathbb{Z}_3$. The linear system $\mathcal{C}_{3,4}=(P_3,\mathcal{L}_3)$ has as set of points to  $P_3=\{(0,1),(1,1),(2,1),(0,2),(1,2),(2,2)\}\cup\{p\}\cup\{q\}$ and as set of lines to $\mathcal{L}_3=\mathcal{L}\cup\mathcal{L}_p\cup\mathcal{L}_q$, where
	\begin{eqnarray*}
		\mathcal{L}&=&\{\{(0,1),(1,1),(2,1)\},\{(0,2),(1,2),(2,2)\}\},\\
		\mathcal{L}_p&=&\{\{(0,1),(0,2),p\},\{(1,1),(1,2),p\},\{(2,1),(2,2),p\}\},\\
		\mathcal{L}_q&=&\{\{(1,1),(2,2),q\},\{(0,1),(1,2),q\},\{(0,2),(2,1),q\}\}\end{eqnarray*}
	and depicted in Figure \ref{fig:nu_2_3}. This linear system is isomorphic to the linear system giving in \cite{MR3727901} Figure 3, which is the linear system with the less number of lines and maximum degree 3 such that $\tau=\nu_2=4$.
\end{example}
\begin{proposition}\label{prop:transversaltau=nu2}
	The linear system $\mathcal{C}_{n,n+1}$ satisfies that: 
	$$\tau(\mathcal{C}_{n,n+1})=n+1$$
\end{proposition}
\paragraph{Proof}
Notice that $\tau(\mathcal{C}_{n,n+1})\leq n+1$ since $\{x_g:\mbox{ $x_g$ is any point of $L_g\in\mathcal{L}$}\}\cup\{p,q\}$ is a transversal of $\mathcal{C}_{n,n+1}$. To prove that $\tau(P_n,\mathcal{L}_n)\geq n+1$, suppose on the
contrary that $\tau(P_n,\mathcal{L}_n)=n$. If $T$ is a transversal of cardinality $n$ then $T\subseteq\Gamma\times	\Gamma\setminus\{e\}$, i.e., $p,q\not\in T$ because, in other case, if $p\in T$ then, by the Pigeonhole principle, there is a line $l_{q_a}\in\mathcal{L}_q$ such that $T\cap l_{q_a}=\emptyset$, since $deg(q)=n$, which is a contradiction, unless that $q\in T$, which implies that there exists $L\in\mathcal{L}$ such that $L\cap T=\emptyset$ (because $|\mathcal{L}|=n-1$), which is also a contradiction. Therefore $T\subseteq\Gamma\times\Gamma\setminus\{e\}$.

Suppose that:
$$T=\{(h_0,f_{g_0}(h_0)),\ldots,(h_{n-1},f_{g_{n-1}}(h_{n-1}))\},$$ where $\{h_0,\ldots,h_{n-1}\}=\{g_0,\ldots,g_{n-1}\}=\Gamma$ and $f_{g_i}=h_i+g_i\neq e$, for $i=0,\ldots,n-1$. Then:
$$\sum_{i=0}^{n-1}f_{h_i}(g_i)=\sum_{i=0}^{n-1}(g_i+h_i)=\sum_{i=0}^{n-1}g_i+\sum_{i=0}^{n-1}h_i=e,$$
since $\displaystyle\sum_{g\in\Gamma}g=\displaystyle\sum_{g\in \Gamma\setminus\{e\}}g=e$, which implies that there exists $f_{h_j}(g_j)\in T$ that satisfies $f_{h_j}(g_j)=e$, which is a contradiction, and consequently $\tau(\mathcal{C}_{n,n+1})=n+1$.\qed

\begin{proposition}\label{prop:2acoplamiento=tau}
	The linear system $\mathcal{C}_{n,n+1}$ satisfies that: 
	$$\nu_2(\mathcal{C}_{n,n+1})=n+1$$.
\end{proposition}
\paragraph{Proof}
Notice that $\nu_2(\mathcal{C}_{n,n+1})\geq n+1$ because, for any two lines $l_{p_g},l_{p_h}\in\mathcal{L}_p$, $\mathcal{L}\cup\{l_{p_g},l_{p_h}\}$ is a 2-packing. To prove that $\nu_2(\mathcal{C}_{n,n+1})\leq n+1$, suppose on the
contrary that $\nu_2(\mathcal{C}_{n,n+1})=n+2$, and that $R$ is a maximum 2-packing of size $n+2$, we analyze to cases:

\textbf{Case $(i)$:} Suppose that
$R=\mathcal{L}\cup\{l_{p_a},l_{p_b},l_{q_c}\}$, where $l_{p_a},l_{p_b}\in\mathcal{L}_p$ and $l_{q_c}\in\mathcal{L}_q$; since there is an unique line $l_p\in\mathcal{L}_p$ which intersect to $l_{q_c}$, then we assume that $l_{p_a}\cap l_{q_c}\neq\emptyset$. By construction of $\mathcal{C}_{n,n+1}$ there exits $L\in\mathcal{L}$ that satisfies
$l_{p_a}\cap l_{q_c}\cap L\neq\emptyset$, inducing a triple point, which is a contradiction.

\textbf{Case $(ii)$:} Let $k$ be an element of $\Gamma\setminus\{e\}$ and 
$R=\{l_{p_a},l_{p_b},l_{q_c},l_{q_d}\}\cup\mathcal{L}\setminus\{L_k\}$ with $l_{p_a},l_{p_b}\in\mathcal{L}_p$ and $l_{q_c},l_{q_d}\in\mathcal{L}_q$, without loss of generality,  suppose that $l_{p_a}\cap l_{q_c}\neq\emptyset$, $l_{p_b}\cap l_{q_d}\neq\emptyset$, $l_{p_a}\cap l_{q_d}=\emptyset$ and
$l_{p_b}\cap l_{q_c}=\emptyset$, otherwise, $R$ is not a 2-packing. It is claimed that there exists $L\in\mathcal{L}\setminus\{L_k\}$ such that either $l_{p_a}\cap l_{q_c}\cap L\neq\emptyset$ or $l_{p_b}\cap l_{q_d}\cap L\neq\emptyset$, which implies that $R$ induce a triple point, which is contradiction and hence $\nu_2(\mathcal{C}_{n,n+1})=n+1$. To verify the claim suppose on the contrary that every 
$L\in\mathcal{L}\setminus\{L_k\}$ satisfies $l_{p_a}\cap l_{q_c}\cap
L=\emptyset$ and $l_{p_b}\cap l_{q_d}\cap L=\emptyset$. It means that $l_{p_a}\cap l_{q_c}\cap L_k\neq\emptyset$ and $l_{p_b}\cap
l_{q_d}\cap L_k\neq\emptyset$. By construction of $\mathcal{C}_{n,n+1}$ it follows that:
\begin{eqnarray*}
	l_{p_i}&=&\{(i,x): x\in\Gamma\setminus\{e\}\}, \mbox{ for all } i\in\Gamma,\\
	l_{q_j}&=&\{(x,x+j): x\in\Gamma\setminus\{e\} \mbox{ and } x+j\neq e\}, \mbox{ for all } j\in\Gamma,\mbox{ and}\\
	L_k&=&\{(x,k): x\in\Gamma\}.
\end{eqnarray*}

If $l_{p_a}\cap l_{q_c}\cap
L_k\neq\emptyset$ and $l_{p_b}\cap l_{q_d}\cap L_k\neq\emptyset$,
then $a+c=b+d=k$. On the other hand, as $l_{p_a}\cap
l_{q_d}=\emptyset$ and $l_{p_b}\cap l_{q_c}=\emptyset$, then
$a+d=b+c=e$. As a consequence of $a+c=b+d=k$ and $a+d=b+c=e$ we obtain
$2k=e$, which is a contradiction. Therefore, $\nu_2(\mathcal{C}_{n,n+1})=n+1$.\qed

Hence, by Proposition \ref{prop:transversaltau=nu2} and Proposition \ref{prop:2acoplamiento=tau} it was proved that:

\begin{theorem}\label{thm:familiacucaracha}
	Let $n=2k+1$, with $k\in\mathbb{N}$, then $$\tau(\mathcal{C}_{n,n+1})=\nu_2(\mathcal{C}_{n,n+1})=n+1,$$with smallest possible cardinality of lines.
\end{theorem}

\subsection{Straight line systems}
A \emph{straight line representation} on $\mathbb{R}^{2}$ of a
linear system $(P,\mathcal{L})$ maps each point $x\in P$ to a point
$p(x)$ of $\mathbb{R}^{2}$, and each line $L\in\mathcal{L}$ to a
straight line segment $l(L)$ of $\mathbb{R}^{2}$ in such a way that
for each point $x\in P$ and line $L\in\mathcal{L}$ satisfies $p(x)\in l(L)$ if and only if $x\in L$, and for each pair of distinct lines	$L,L^\prime\in\mathcal{L}$ satisfies $l(L)\cap l(L^\prime)=\{p(x):x\in L\cap L^\prime\}$. A
\emph{straight line system} $(P,\mathcal{L})$ is a linear system,
such that it has a straight line representation on $\mathbb{R}^{2}$. In \cite{MR3727901} was proved that the linear system $\mathcal{C}_{3,4}$ is not a straight one. The \emph{Levi graph} of a linear system $(P,\mathcal{L})$, denoted by $B(P,\mathcal{L})$, is a bipartite graph with vertex set
$V=P\cup\mathcal{L}$, where two vertices $p\in P$, and $L\in\mathcal{L}$ are adjacent if and only if $p\in L$.

In the same way as in \cite{MR3727901} and according to \cite{Kaufmann:2009:SDH:1506879.1506920}, any straight line system is Zykov-planar, see also \cite{MR0401556}. Zykov proposed to represent the lines of a set system by a subset of the faces of a planar map on $R^2$, i.e., a set system $(X,\mathcal{F})$ is Zykov-planar if there exists a planar graph $G$ (not necessarily a simple graph) such that $V(G)=X$ and $G$ can be drawn in the plane with faces of $G$ two-colored (say red and blue) so that there exists a bijection between the red faces of $G$ and the subsets of
$\mathcal{F}$ such that a point $x$ is incident with a red face if and only if it is incident with the corresponding subset. In \cite{MR0360328} was shown that the Zykov's definition is equivalent to the following: A set system $(X,\mathcal{F})$ is Zykov-planar if and only if the Levi graph $B(X,\mathcal{F})$ is planar. It is well-known that for any planar graph $G$ the size of $G$, $|E(G)|$,  is upper bounded by $\frac{k(|V(G)|-2)}{k-2}$ (see \cite{Bondy} page 135, exercise 9.3.1 (a)), where $k$ is the \emph{girth} of $G$ (the length of a shortest cycle contained in the graph $G$). It is not difficult to prove that the Levi graph $B(\mathcal{C}_{n,n+1})$ of $\mathcal{C}_{n,n+1}$ is not a planar graph, since the size of the girth of $B(\mathcal{C}_{n,n+1})$ is $6$, it follows: $$3n^2-n=|E(\mathcal{C}_{n,n+1})|>\frac{3(n^2+2n-1)}{2},$$for all $n\geq3$. 
Therefore, the linear system $\mathcal{C}_{n,n+1}$ is not a straight line system. 

Finally, as a Corollary of Theorem \ref{thm:casicucarachageneral}, we have the following: 

\begin{corollary}
	Let $(P,\mathcal{L})$ be a straight line system with $p,q\in P$ be two points such that $deg(p)=\Delta(P,\mathcal{L})$ and $deg(q)=\max\{deg(x): x\in
	P\setminus\{p\}\}$. If $|\mathcal{L}|\leq deg(p)+deg(q)+\nu_2(P,\mathcal{L})-3$, then $\tau(P,\mathcal{L})\leq\nu_2(P,\mathcal{L})-1$.
\end{corollary}
\section{Intersecting $r$-uniform linear systems with $\tau=\nu_2=r$}\label{sec:intersectante_igualad}

In this subsection, we give some properties of $r$-uniform linear systems that satisfies $\tau=\nu_2=r$ as well as a characterization of $4$-uniform linear systems with $\tau=\nu_2=4$.

Let $\mathbb{L}_r$ be the family of intersecting linear systems $(P,\mathcal{L})$ of rank $r$ that satisfies $\tau(P,\mathcal{L})=\nu_2(P,\mathcal{L})=r$, then we have the following lemma:

\begin{lemma}
	Each element of $\mathbb{L}_r$ is an $r$-uniform linear system. 
\end{lemma}
\paragraph{Proof}
Let consider $(P,\mathcal{L})\in\mathbb{L}_r$ and $l\in\mathcal{L}$ any line of $(P,\mathcal{L})$. It is clear that $T=\{p\in l: deg(p)\geq2\}$ is a transversal of $(P,\mathcal{L})$. Hence $r=\tau(P,\mathcal{L})\leq|T|\leq r$, which implies that $|l|=r$, for all $l\in\mathcal{L}$. Moreover, $deg(p)\geq2$, for all $p\in l$ and $l\in\mathcal{L}$. \qed

In \cite{DHLK} was proved the following:

\begin{lemma}\cite{DHLK}
	Let $(P,\mathcal{L})$ be an $r$-uniform intersecting linear system then every edge of $(P,\mathcal{L})$ has at most one vertex of degree 2. Moreover $\Delta(P,\mathcal{L})\leq r$. 
\end{lemma}

\begin{lemma}\label{lemma:tamano_lineas}\cite{DHLK}
	Let $(P,\mathcal{L})$ be an $r$-uniform intersecting linear system then $$3(r-1)\leq|\mathcal{L}|\leq r^2-r+1.$$
\end{lemma}

Hence, by Theorem \ref{thm:casicucarachageneral} and Lemma \ref{lemma:tamano_lineas} it follows:

\begin{corollary}\label{coro:nuevo}
	If $(P,\mathcal{L})\in\mathbb{L}_r$ then $3(r-1)+1\leq|\mathcal{L}|\leq r^2-r+1$.
\end{corollary}

In \cite{MR3727901} was proved that the linear systems $(P,\mathcal{L})$ with $|\mathcal{L}|>\nu_2(P,\mathcal{L})$ and $\nu_2(P,\mathcal{L})\in\{2, 3, 4\}$ satisfy $\tau(P,\mathcal{L})\leq\nu_2(P,\mathcal{L})$; and when attain the equality, they are a special family of linear subsystems of the projective plane of order $3$, $\Pi_3$ (some of them 4-uniform intersecting linear systems) with transversal and $2$-packing numbers equal to $4$. 
Recall that a \emph{finite projective plane} (or merely \emph{projective plane}) is a linear system satisfying that any pair of points have a common line, any pair of lines have a common point and there exist four points in general
position (there are not three collinear points). It is well known that, if $(P,\mathcal{L})$ is a projective plane, there exists a number $q\in\mathbb{N}$, called \emph{order of projective plane}, such that every point (line, respectively) of $(P,\mathcal{L})$ is incident to exactly $q+1$ lines (points, respectively), and $(P,\mathcal{L})$ contains exactly $q^2+q+1$ points (lines, respectively). In addition to this, it is well known that projective planes of order $q$, denoted by $\Pi_q$, exist when $q$ is a power prime. For more information about the existence and the unicity of projective planes see, for instance, \cite{B86,B95}.

Given a linear system $(P,\mathcal{L})$, a
\emph{triangle} $\mathcal{T}$ of $(P,\mathcal{L})$, is the linear
subsystem of $(P,\mathcal{L})$ induced by three points in general position (non collinear) and the three lines induced by them. In \cite{MR3727901} was defined $\mathcal{C}=(P_\mathcal{C},\mathcal{L}_\mathcal{C})$ to be the linear system obtained from $\Pi_3$ by deleting $\mathcal{T}$; also there was defined $\mathcal{C}_{4,4}$ to be the family of linear systems $(P,\mathcal{L})$ with $\nu_2(P,\mathcal{L})=4$, such that:
\begin{itemize}
	\item[i)] $\mathcal{C}$ is a linear subsystem of $(P,\mathcal{L})$; and
	
	\item[ii)] $(P,\mathcal{L})$ is a linear subsystem of $\Pi_3$,
\end{itemize}

this is $\mathcal{C}_{4,4}=\{(P,\mathcal{L}):\mathcal{C}\subseteq(P,\mathcal{%
	L})\subseteq\Pi_3\mbox{ and $\nu_2(P,\mathcal{L})=4$}\}$.

Hence, the authors proved the following:

\begin{theorem}\cite{MR3727901}\label{teo:anterior}
	Let $(P,\mathcal{L})$ be a linear system with $\nu_2(P,\mathcal{L})=4$. Then, $\tau(P,\mathcal{L})=\nu_2(P,\mathcal{L})=4$ if and only if $(P,\mathcal{L})\in\mathcal{C}_{4,4}$. 	
\end{theorem}
Now, consider the projective plane
$\Pi_3$ and a triangle $\mathcal{T}$ of $\Pi_3$ (see $(a)$ of Figure \ref{fig:planos}). Define $\mathcal{\hat{C}}=(P_\mathcal{C},\mathcal{L}_\mathcal{C})$ to be the linear subsystem induced by $\mathcal{L}_\mathcal{C}=\mathcal{L}\setminus\mathcal{T}$ (see $(b)$ of Figure \ref{fig:planos}). The linear system $\mathcal{\hat{C}}=(P_\mathcal{C},\mathcal{L}_\mathcal{C})$ just defined has ten points and ten lines. Define $\mathcal{\hat{C}}_{4,4}$ to be the family of 4-uniform intersecting linear systems $(P,\mathcal{L})$ with $\nu_2(P,\mathcal{L})=4$, such that:

\begin{figure}
	\begin{center}
		\subfigure[]{\includegraphics[height =3.5cm]{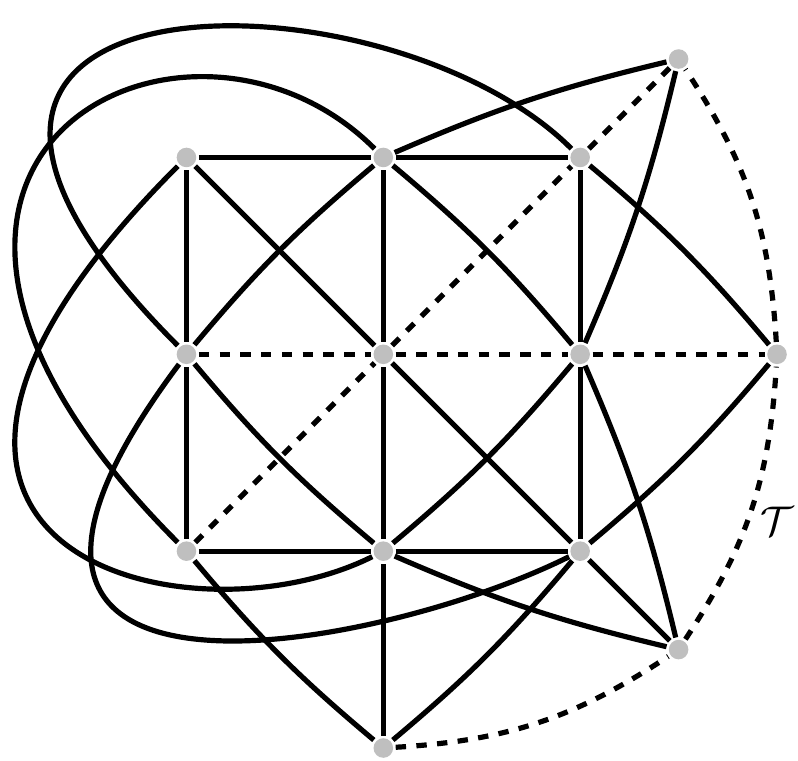}}
		\hspace{1cm}
		\subfigure[]{\includegraphics[height =3.5cm]{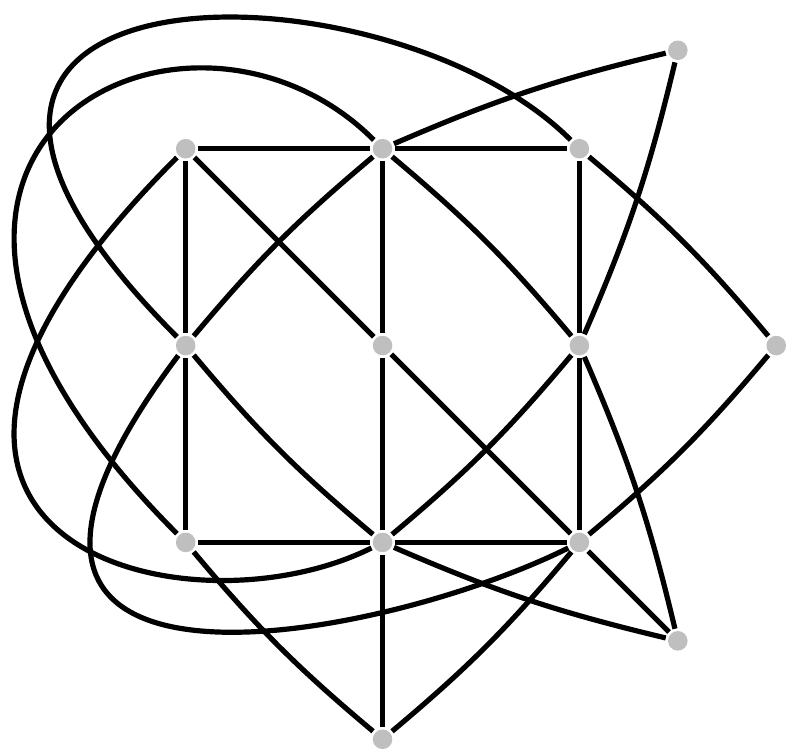}}
	\end{center}
	\caption{(a) Projective plane of order 3, $\Pi_3$ and $(b)$ Linear system obtained from $\Pi_3$ by deleting the lines of the triangle $\mathcal{T}$.}\label{fig:planos}
\end{figure}

\begin{itemize}
	\item[i)] $\mathcal{\hat{C}}$ is a linear subsystem of $(P,\mathcal{L})$; and
	
	\item[ii)] $(P,\mathcal{L})$ is a linear subsystem of $\Pi_3$,
\end{itemize}

It is clear that $\mathcal{\hat{C}}_{4,4}\subseteq\mathcal{C}_{4,4}$ and each linear system $(P,\mathcal{L})\in\mathcal{\hat{C}}_{4,4}$ is an 4-uniform intersecting linear system. Hence

\begin{corollary}
	$(P,\mathcal{L})\in\mathbb{L}_4$ if and only if $(P,\mathcal{L})\in\mathcal{\hat{C}}_{4,4}$.
\end{corollary}

\end{document}